\documentclass[12pt]{article}
\usepackage{}
\usepackage{amsfonts}

\usepackage{amsmath, amssymb}
\usepackage{amssymb}
\def\mineappendix{
        \setcounter{section}{1}
        \setcounter{subsection}{0}
        \def\thesection{\Alph{section}}
        \def\sectionap{\@startsection  {section}{1}{\z@}
                        {-3.5ex plus-1ex minus-.2ex} {0ex plus.2ex}
                        {\reset@font\Large\bf  Appendix:  \, }
                        }
        }
\makeatother
\def\Proclaim #1. #2\par{\bigbreak\noindent{\sc#1.\enspace}{\it#2}\par}

\newcommand{\gwi}[1]{\left< \, #1 \, \right>_{0}}
\newcommand{\gwid}[1]{\left< \, #1 \, \right>_{0,d}}

\newcommand{\gwii}[1]{\left< \hspace{-2pt} \left< \, #1 \,
        \right>  \hspace{-2pt} \right>_{0}}

\newcommand{\gwione}[1]{\left< \, #1 \, \right>_{1}}
\newcommand{\gwioned}[1]{\left< \, #1 \, \right>_{1,d}}
\newcommand{\gwioneo}[1]{\left< \, #1 \, \right>_{1,0}}
\newcommand{\gwiione}[1]{\left< \hspace{-2pt} \left< \, #1 \,
        \right> \hspace{-2pt} \right>_{1}}

\newcommand{\gwitwo}[1]{\left< \, #1 \, \right>_{2}}
\newcommand{\gwitwod}[1]{\left< \, #1 \, \right>_{2,d}}

\newcommand{\gwig}[1]{\left< \, #1 \, \right>_{g}}
\newcommand{\gwigd}[1]{\left< \, #1 \, \right>_{g,d}}

\newcommand{\gwiig}[1]{\left< \hspace{-2pt} \left< \, #1 \,
    \right> \hspace{-2pt} \right>_{g}}

\newcommand{\grav}[2]{\tau_{#1}(\gamma_{#2})}
\newcommand{\grava}[1]{\tau_{#1}(\gamma_{\alpha})}

\newcommand{\ga}{\gamma_{\alpha}}
\newcommand{\gua}{\gamma^{\alpha}}
\newcommand{\gb}{\gamma_{\beta}}
\newcommand{\gub}{\gamma^{\beta}}

\newcommand{\vs}{{\cal S}}

\newcommand{\ve}{{\mathcal E}}

\newcommand{\qp}{\circ}

\newtheorem{lem}{Lemma}[section]
\newtheorem{cor}[lem]{Corollary}
\newtheorem{thm}[lem]{Theorem}

\newtheorem{rem}[lem]{Remark}

\title{ Genus-2 G-function for $\mathbb{P}^{1}$
orbifolds}

\author{Xiaobo Liu \thanks{Research of the first author was partially supported by
NSFC Tianyuan special fund 11226027, SRFDP grant
20120001110051, and Peking University 985 fund.}, \quad  Xin Wang
\thanks{Research of the second author was partially supported by
SRFDP grant 20120001110051.}}

\date{}

\begin{document}
\maketitle
\begin{abstract}
In this paper we prove that for Gromov-Witten theory of
$\mathbb{P}^{1}$ orbifolds of ADE type the genus-2 G-function
introduced by B. Dubrovin, S. Liu, and Y. Zhang vanishes. Together
with our results in \cite{LW}, this completely solves the main
conjecture in their paper \cite{DLZ}. In the process, we also found
a sufficient condition for the vanishing of the genus-2 G-function
which is weaker than the condition given in our previous paper
\cite{LW}.

\end{abstract}

The genus-2 G-function for semisimple Frobenius manifolds was
introduced by B. Dubrovin, S. Liu, and Y. Zhang in \cite{DLZ}. On a semisimple
Frobenius manifold, there are two coordinate systems: One is called the
flat coordinate which is characterized by the property that the non-degenerate
pairing on the tangent bundle can be represented by a constant matrix
using this coordinate; Another one is called the canonical coordinate
in which the coordinate vector fields are idempotents of the product structure.
Let $F_2$ be the genus-2 potential function for a semisimple Frobenius
manifold. In \cite{DLZ},
$F_{2}$ was written as a summation of two functions:
 \begin{align}
 F_{2}=F'_{2}+G^{(2)}.    \label{F2}
 \end{align}
In this decomposition, the simple part
$F'_{2}$ is a function which only depends on genus-0 and genus-1 data and
can be expressed purely in terms of flat coordinates on the
big phase space of Frobenius manifolds. We do not need $F'_{2}$ in this paper.
Precise definition of $F'_2$ can be found in \cite{DLZ} where a graphical
representation of this function was also given.
The function $G^{(2)}$ appeared in decomposition \eqref{F2} is called the {\bf genus-$2$ G-function}.
This is a very complicated function
which depends on the canonical coordinate system of semisimple Frobenius manifolds.
The precise definition of $G^{(2)}$ can be found in Appendix~\ref{sec:G2}.

For each cohomological field theory, genus-0 part of the theory defines a Frobenius
manifold (cf. \cite{KM}). In this case $F_2$ is the
generating function for genus-2 descendant invariants. These
geometric invariants appear as coefficients in the formal power series representation
of $F_2$ in terms of flat coordinates. In general, explicit formulas for the transition functions
between flat coordinate and canonical coordinate for semisimple Frobenius manifolds
are not easy to obtain. Therefore for the computation of geometric invariants, it is desirable
to have a formula for $F_2$ only using flat coordinates. If the genus-2 G-function is $0$,
this can be easily achieved due to decomposition \eqref{F2}. In \cite{DLZ},
it was conjectured that the genus-2 G-function vanishes for the cohomological field theories
associated to ADE singularities and for Gromov-Witten theory of $\mathbb{P}^{1}$ orbifolds of ADE type.
For $A_n$ singularities, a proof of this conjecture was given by Y. Fu, S. Liu, Y. Zhang, and
C. Zhou in \cite{FLZZ}, which relies heavily on the specific properties of the Frobenius manifold
structure of $A_n$ singularities. A more geometric proof of this conjecture for all ADE singularities
was given in our previous paper \cite{LW}. In this paper we prove this conjecture for
$\mathbb{P}^{1}$ orbifolds of ADE type, and thus completely solve the main conjecture in
\cite{DLZ}. More precisely, we have the following
\begin{thm} \label{thm:ADEG}
The genus-2 G-function vanishes for the Gromov-Witten theory of $\mathbb{P}^{1}$ orbifolds of $ADE$ type.
\end{thm}

In the process of proving this theorem, we also found a sufficient condition for the vanishing of
the genus-2 G-function which is weaker than the condition given in Theorem 0.1 of \cite{LW}.
More precisely, each cohomological field theory defines a Frobenius manifold structure on
a vector space $\mathcal{H}$ with a non-degenerate pairing $\eta$.
Fix a basis $\{\gamma_{1},\gamma_{2},...,\gamma_{N}\}$ of $\mathcal{H}$
with $\gamma_1$ being the identity for the Frobenius manifold structure on $\mathcal{H}$.
Let $\eta_{\alpha\beta}$ be entries for the matrix of $\eta$
with respect to this basis. Entries for inverse matrix of $\eta$ are denoted
$\eta^{\alpha\beta}$.
For any $\alpha$, $\gamma^{\alpha}$ is
defined by $\gamma^{\alpha}=\sum_{\beta}\eta^{\alpha\beta}\gamma_{\beta}$. Let
$\gwig{\grav{n_1}{\alpha_1} \cdots \grav{n_k}{\alpha_k}}$ be
the genus-$g$ descendant invariants of the cohomological field theory.
We will identify $\grava{0}$ with $\gamma_{\alpha}$ and call
$\gwig{\gamma_{\alpha_1} \cdots \gamma_{\alpha_k}}$ genus-$g$ primary invariants.
The Frobenius manifold structure on $\mathcal{H}$ is determined by the
genus-$0$ primary invariants. In \cite{LW}, we proved that
the genus-2 G-function vanishes under the assumption that
certain invariants of genus $g \leq 2$ are zero. In this paper, we
consider three similar conditions:
\begin{align*}
&(C1) \quad\quad
\sum_{\alpha,\beta}\gwi{ {\gamma _\alpha }{\gamma ^\alpha }{\gamma _\beta }{\gamma ^\beta }\gamma_{\alpha_{1}}\cdot\cdot\cdot \gamma_{\alpha_{k}}} =0,
\\&(C2) \quad\quad  \gwione{ \gamma_{\alpha_{1}}\cdot\cdot\cdot \gamma_{\alpha_{k}}} =0 \quad(k\geq2) \quad \mbox{and} \quad
\gwione{\gamma_{\alpha}}=C\eta_{1\alpha},
\\&(C3) \quad\quad  \gwitwo{ \gamma_{\alpha_{1}}\cdot\cdot\cdot \gamma_{\alpha_{k}}} =0 \quad \mbox{and} \quad
\gwitwo{ \tau_{1}(\gamma_{\alpha_{1}})\gamma_{\alpha_{2}}\cdot\cdot\cdot \gamma_{\alpha_{k}}} =0,
\end{align*}
for all $\alpha_{1},\cdot\cdot\cdot,\alpha_{k}$, $k\geq 1$, where $C$ is a constant which only
depends on $\mathcal{H}$. Conditions (C1) and (C3) are exactly the
same as the corresponding conditions in Theorem 0.1 of \cite{LW}.
Condition (C2) is weaker than the corresponding condition in \cite{LW} in the sense that
the latter one
is equivalent to setting $C=0$ in condition (C2). In this paper we will also prove the following
\begin{thm} \label{thm:main}
For any semisimple cohomological field theory satisfying conditions (C1),(C2) and (C3),
 the genus-$2$ G-function vanishes.
\end{thm}

To prove Theorem \ref{thm:ADEG}, we need to show that conditions (C1)--(C3) are satisfied
for the Gromov-Witten theory of $\mathbb{P}^{1}$ orbifolds of $ADE$ type.
As a byproduct of this process, we also proves another conjecture, i.e. Conjecture 3.20, in \cite{DLZ}
 up to a constant (cf. Corollary~\ref{cor:P1ADEg=1G}).

This paper is organized as follows.
In section~\ref{sec:prel}, we introduce notations and review basic properties of
idempotents and rotation coefficients needed in this paper.
We prove Theorem~\ref{thm:main} in section~\ref{sec:pfmain}.
In section~\ref{sec:pfcor} we give the proof of Theorem~\ref{thm:ADEG}.
In appendix, we recall the precise definition for the genus-2 G-function.

\section{Preliminaries}
\label{sec:prel}

In this section we set up notations and review some basic structures for semisimple Frobenius
manifold $\mathcal{H}$ and their extensions to the {\it big phase space}
$\mathcal{P}:=\prod_{n=0}^{\infty} \mathcal{H}$. The first copy of $\mathcal{H}$
in the product defining $\mathcal{P}$ is usually identified with $\mathcal{H}$ itself,
and is called the {\it small phase space}.
A basis $\{\gamma_1, \cdots, \gamma_n\}$ of $\mathcal{H}$ induces a basis
$\{\tau_n(\gamma_{\alpha}) \mid 1 \leq \alpha \leq N, \,\,\, n \geq 0\}$ of $\mathcal{P}$.
The coordinates on $\mathcal{P}$ with respected to this basis is denoted
by $(t_{n}^{\alpha} \mid 1 \leq \alpha \leq N, \,\,\, n \geq 0)$. This is called the
{\it flat coordinate system} on $\mathcal{P}$.
Each $\tau_n(\gamma_{\alpha})$ can also be viewed as a coordinate vector field on $\mathcal{P}$
and we identify $\tau_0(\gamma_{\alpha})$ with $\gamma_{\alpha}$. Linear combinations of
$\gamma_{\alpha}$ with coefficients being functions on $\mathcal{P}$ are called {\it primary
vector fields}.
The product structure on Frobenius manifold $\mathcal{H}$ can be naturally extended to
an associative product,
denoted by "$\circ$", on $\mathcal{P}$
using genus-0 potential function $F_0$ on the big phase space (cf. \cite{L02}).
For a cohomological field theory, $F_0$
is the generating function for genus-0 descendent invariants.
The product "$\circ$" can be used to define the quantum product of vector fields on $\mathcal{P}$.
Semisimplicity of $\mathcal{H}$
guarantees the existence of mutually commuting primary vector fields $ \ve_{1}, \cdots, \ve_{N} $ such that
\[ \ve_i \circ \ve_j = \delta_{ij} \,\, \ve_i.\]
These vector fields are called {\it idempotents}.
On the small phase space, $\ve_i$ can be identified with $\frac{\partial}{\partial u_i}$
where $(u_1, \cdots, u_N)$ is the {\it canonical coordinate} system on $\mathcal{H}$
(cf. \cite{D}).
Functions $u_i$ can also be extended to functions on the big phase space $\mathcal{P}$ as the eigenvalue
of the quantum
multiplication by the Euler vector field on $\mathcal{P}$ (cf. \cite{L06}).
These functions satisfy the property
\[ \ve_i u_j = \delta_{ij}. \]
The non-degenerate pairing $\eta$ on $\mathcal{H}$ can also be extended to a bilinear form
$<\cdot, \, \cdot>$ which is only non-degenerate on the space of primary vector fields
(cf. \cite{L06}). Idempotents are orthogonal with respect to $<\cdot, \, \cdot>$. Define
$$g_{i}:=<\mathcal{E}_{i},\mathcal{E}_{i}>,\quad
h_{i}=\sqrt{g_{i}}$$
and set
\begin{align}
r_{ij}:=\frac{\mathcal{E}_{j}}{\sqrt{g_{j}}}\sqrt{g_{i}}. \label{r_{ij}}
\end{align}
Functions $r_{ij}$ are called the {\it rotation coefficients on the big phase space} (cf. \cite{L06}).
These functions are symmetric with respect to $i$ and $j$.
For $i\neq j$, the restriction of $r_{ij}$ to the small phase space coincides with Dubrovin's
 definition of
rotation coefficients for semisimple Frobenius manifolds (cf. \cite{D}). For $i=j$, Dubrovin
set the corresponding rotation coefficients to 0. This is different from our definition
of $r_{ii}$ given by equation \eqref{r_{ij}}.
Except in Appendix \ref{sec:G2}, throughout this paper we will use equation~\eqref{r_{ij}}
for rotation coefficients instead of Dubrovin's definition.
The {\it string equation} on the small phase space $\mathcal{H}$
can be written as
\begin{align}
\sum_{j}r_{ij}h_{j} \mid_{\mathcal{H}}=0 \label{eq:string}
\end{align} for all $i$.

We now review basic properties of rotation coefficients which are most relevant to
our calculations in this paper. More properties about idempotents and rotation coefficients on
the big phase spaces can be found in \cite{L06} and \cite{L07}. Although results in
\cite{L06} and \cite{L07} were stated for Gromov-Witten theory, they can be easily extended
to all semisimple cohomological field theories.

As observed in \cite{LW}, a key step in proving the vanishing of the genus-2 G-function
is to express this function in terms of $r_{ij}$, $h_i$, and the following functions
\begin{align*}
v_{ij}:=(u_{j}-u_{i})r_{ij},
\end{align*}
$$
\theta_{ij}:=\frac{1}{u_{j}-u_{i}}\Big(r_{ij}+\sum_{k}{r_{ik}v_{jk}}
\Big)
$$
and
$$\Omega_{ij}:
=\frac{1}{u_{j}-u_{i}}\Big(\theta_{ij}-\theta_{ji}
+\sum_{k,l}r_{il}r_{jk}v_{kl}\Big).$$
$\theta_{ij}$ and $\Omega_{ij}$ are only defined for $i\neq j$.
Obviously,
\begin{align}\theta_{ij}+\theta_{ji}=-\sum_{k}r_{ik}r_{jk},
\hspace{20pt} \Omega_{ij}=\Omega_{ji}
\label{eq:theta}\end{align} for any $i\neq j$.
We might consider $\theta_{ij}$ and $\Omega_{ij}$ as functions having poles of
order 1 and 2 respectively in terms of $u_{1},...,u_{N}$.
These functions appear naturally when taking derivatives of $r_{ij}$ along
idempotents. More precisely
\begin{equation}
\mathcal{E}_{k}r_{ij}=r_{ik}r_{jk}+
\begin{cases}
0& \text{if $i\neq j\neq k$},\\
\theta_{ij}& \text{if $k=i\neq j$},\\
\sqrt{\frac{g_{k}}{g_{i}}}\theta_{ik}& \text{if $i=j\neq k$},\\
-2\sum_{l}{r_{il}^{2}}+\sum_{p\neq i}{\sqrt{\frac{g_{p}}{g_{i}}}\theta_{pi}}
+\frac{1}{g_{i}}<\tau_{-}^{2}(\mathcal{S}),\mathcal{E}_{i}>,& \text{if $i=j=k$}  \label{eq:der of r}
\end{cases}
\end{equation}
and
$$\mathcal{E}_{j}\theta_{ij}=\Big(r_{jj}-\sqrt{\frac{g_{j}}{g_{i}}}r_{ij}\Big)
\theta_{ij}-\Omega_{ij}$$ for $i\neq j$.
The vector field $\vs$ in equation \eqref{eq:der of r} is the {\it string vector field}
defined by
\[ \vs := \gamma_1 - \sum_{n, \alpha} t_{n}^{\alpha} \tau_{n-1}(\gamma_{\alpha}). \]
The operator $\tau_{-}$ acting on the space of vector fields  is given by
$\tau_{-}(\tau_n(\gamma_{\alpha})) = \tau_{n-1}(\gamma_{\alpha})$.
Notice that on the small phase space $\tau_{-}^{2}(\mathcal{S})=0$.

Conditions (C1)--C(3) can be better represented using
tensors $\gwiig{\cdot\cdot\cdot}$ defined by
\begin{align}
\gwiig{\mathcal{W}_{1}\mathcal{W}_{2}\cdot\cdot\cdot\mathcal{W}_{k}}
:=\sum_{m_{1},\alpha_{1},...,m_{k},\alpha_{k}}f_{m_{1},\alpha_{1}}^{1}
\cdot\cdot\cdot f_{m_{k},\alpha_{k}}^{k}\frac{\partial^{k}}{\partial t_{m_{1}}^{\alpha_{1}}\cdot\cdot\cdot \partial t_{m_{k}}^{\alpha_{k}}}F_{g}
\end{align}
for vector fields
$\mathcal{W}_{i}=\sum_{m,\alpha}f_{m,\alpha}^{i}\frac{\partial}{\partial t_{m}^{\alpha}}$
where $f_{m,\alpha}^{i}$
are functions on the big phase space. This tensor is called $k$-$point$ $(correlation)$ $function$.
As in \cite{L07}, we will use the following notation for genus-1 $k$-point functions:
$$\phi_{i_{1},\cdot\cdot\cdot,i_{k}}:=
\gwiione{\mathcal{E}_{i_{1}},\cdot\cdot\cdot,\mathcal{E}_{i_{k}}}.$$
It was proved in \cite{L06} that the genus-1 1-point functions are given by
\begin{align} 24\phi_{i}=-12\sum_{j}{r_{ij}v_{ij}}
-\sum_{j}{\frac{h_{i}}{h_{j}}}r_{ij}
\label{eq:g=1 1-pt}
\end{align} for all $i$.
Higher point genus-1 functions can also be computed using recursion formulas given in \cite{L06}.
For example, the genus-1 2-point functions can be written as
\begin{align}
24\phi_{ij}=& 12r_{ij}^2 + \sum\limits_l {({r_{il}}{r_{jl}}{{{h_i}{h_j}} \over {h_l^2}} - {r_{ij}}{r_{il}}{{{h_j}} \over {{h_l}}} - {r_{ij}}{r_{jl}}{{{h_i}} \over {{h_l}}})}-\{{\theta _{ij}}{{{h_j}} \over {{h_i}}} + {\theta _{ji}}{{{h_i}} \over {{h_j}}}\}
\nonumber\\&-24r_{ij}\{\frac{h_{j}}{h_{i}}\phi_{i}+\frac{h_{i}}{h_{j}}\phi_{j}\}
\label{eq:g=1 2-ptij}
\end{align}
for $i\neq j$.
Explicit formulas for some genus-1 3-point functions can also be found
in \cite{LW}.

\section{Proof of Theorem \ref{thm:main}}
\label{sec:pfmain}
\allowdisplaybreaks

The only difference between Theorem \ref{thm:main} in this paper and Theorem 0.1 in \cite{LW}
is the appearance of the constant $C$ in condition (C2).  In \cite{LW} this constant $C$ is 0 while
for Theorem~\ref{thm:main} in this paper $C$ may not be zero.
In this section, we will use the same method
as in \cite{LW} to prove Theorem \ref{thm:main}. All formulas in the proof of Theorem 0.1 in \cite{LW}
will be adapted to accommodate this change of condition (C2).
In this process, extra terms containing the constant $C$ will appear in these formulas, and
 terms not containing $C$ are exactly the same as in \cite{LW}. The proof in \cite{LW}
shows that the summation of the contributions  of all terms not containing $C$ to the genus-2 G-function
is $0$. Therefore to prove Theorem \ref{thm:main}, we only need to show that the contributions
of all terms containing $C$ to the genus-2 G-function add up to $0$.

We first note that for semisimple cohomological field theories,
condition (C2) have the following equivalent form which is easier to use:
\begin{lem} \label{lem:C2Idemp}
If $\mathcal{H}$ is semisimple, condition (C2) is equivalent to
\begin{align}
\phi_{i}\mid_{\mathcal{H}}=Ch_{i}^{2} \label{eq:g=1 k=1-pt}
\end{align}
for all $i$ and
\begin{align}
\phi_{i_{1}\cdot\cdot\cdot i_{k}}\mid_{\mathcal{H}}=0 \label{eq:g=1 k-pt}
\end{align}
for all
$i_{1},\cdot\cdot\cdot,i_{k}$ with $k\geq 2$.
\end{lem}
{\bf Proof}:
Since $\eta_{\alpha \beta}$ is constant in flat coordinates, condition $(C2)$ is equivalent to
\begin{equation} \label{eqn:1ptg1flat}
\gwiione{\ga}|_{\mathcal{H}}=C \eta_{1 \alpha}
\end{equation}
for all $\alpha$.
By the string equation
$\eta_{1 \alpha} = \gwii{\gamma_{1}\gamma_{1}\ga}|_{\mathcal{H}}. $
So equation \eqref{eqn:1ptg1flat} can be rewritten as
\[\gwiione{\ga}|_{\mathcal{H}}=C\gwii{\gamma_{1}\gamma_{1}\ga}|_{\mathcal{H}}\]
for all $\alpha$. Since both sides of this equation are linear with respect to $\gamma_{\alpha}$,
we have
\[ \phi_i \mid_{\mathcal{H}} := \gwiione{\mathcal {E}_{i}}|_{\mathcal{H}}
=C\gwii{\gamma_{1}\gamma_{1}\mathcal {E}_{i}}|_{\mathcal{H}}. \]
Since $\gamma_1$ is the identity element of the quantum product on the small phase space,
by definition of idempotents and the associativity of the quantum product
\[ \gwii{\gamma_{1}\gamma_{1}\mathcal {E}_{i}}|_{\mathcal{H}}
= \gwii{\gamma_{1}\gamma_{1} \{ \ve_{i} \qp \ve_i \}} |_{\mathcal{H}}
= \gwii{\gamma_{1} \{ \gamma_{1} \qp \ve_{i} \} \ve_i } |_{\mathcal{H}}
= \gwii{\gamma_{1}\mathcal {E}_{i}\mathcal {E}_{i}}|_{\mathcal{H}}= h_{i}^{2}.\]
In the last equality, we have also used the string equation.
In the semisimple case,
idempotents $\ve_1, \cdots, \ve_N$ form a basis for the space of primary vector fields.
Therefore equation \eqref{eqn:1ptg1flat} is equivalent to
equation \eqref{eq:g=1 k=1-pt}, and both of them are equivalent to condition (C2).

Equation \eqref{eq:g=1 k-pt} is a consequence of either equation \eqref{eqn:1ptg1flat} or
equation \eqref{eq:g=1 k=1-pt}. In fact by equation \eqref{eqn:1ptg1flat},
derivatives of $\gwiione{\ga}|_{\mathcal{H}}$ are all 0.
By linearity of the $k$-tensor $\gwiione{\cdots}$, this implies
\begin{equation} \label{eqn: g1k>1}
\gwiione{w_1 \cdots w_k} |_{\mathcal{H}} = 0
\end{equation}
for any primary vector fields $w_1, \cdots, w_k$ with $k \geq 2$.
This implies equation \eqref{eq:g=1 k-pt}. The lemma is thus proved.
$\Box$

An immediate consequence of this lemma is the following
\begin{cor} Under the condition $(C2)$, on the small phase space the following identities hold:
\begin{align}
&\sum_{j}{r_{ij}v_{ij}} \,\, = \,\,
-\frac{1}{12}\sum_{j}r_{ij}\frac{h_{i}}{h_{j}}-2Ch_{i}^{2}
\label{eq:g=1 1-pt small}
\end{align}
for all $i$,
and
\begin{align}
&{\theta _{ij}}{{{h_j}} \over {{h_i}}} + {\theta _{ji}}{{{h_i}} \over {{h_j}}} \,\, = \,\, 12r_{ij}^2 + \sum\limits_l {({r_{il}}{r_{jl}}{{{h_i}{h_j}} \over {h_l^2}} - {r_{ij}}{r_{il}}{{{h_j}} \over {{h_l}}} - {r_{ij}}{r_{jl}}{{{h_i}} \over {{h_l}}})}-48Cr_{ij}h_{i}h_{j}\label{eq:g=1 2-pt small}
  \end{align}
  for all $i\neq j$.
\end{cor}
{\bf Proof}:
Equation $(\ref{eq:g=1 1-pt small})$ follows from equations (\ref{eq:g=1 1-pt}) and (\ref{eq:g=1 k=1-pt}).
Equation (\ref{eq:g=1 2-pt small}) follows from equations (\ref{eq:g=1 2-ptij}),
(\ref{eq:g=1 k=1-pt}) and
$\phi_{ij}\mid_{\mathcal{H}}=0$ for  $i\neq j$.
$\Box$

Recall that the genus-2 G-function  $G^{(2)}(u,u_{x},u_{xx})$ defined in \cite{DLZ} can be written as:
                  \begin{align}
                   G^{(2)}(u,u_{x},u_{xx})
                   =&\sum_{i}G_{i}^{(2)}(u,u_{x})u_{xx}^{i}
                   +\sum_{i\neq j}G_{ij}^{(2)}(u)\frac{(u_{x}^{j})^{3}}{u_{x}^{i}}
                   \nonumber\\&+\frac{1}{2}\sum_{i,j}P_{ij}^{(2)}(u)u_{x}^{i}u_{x}^{j}
                   +\sum_{i}Q_{i}^{(2)}(u)(u_{x}^{i})^{2}.
                              \end{align}
In this formula, $u=(u_1, \cdots, u_N)$ is the canonical coordinate system on the small phase space
$\mathcal{H}$. The flat coordinate system on $\mathcal{H}$ are given by $t_0^{\alpha} \mid_{\mathcal{H}}$ for
$\alpha=1, \cdots, N$.
Functions on $\mathcal{H}$ can be lifted to functions on the big phase space $\mathcal{P}$
via the transformation
\begin{equation} \label{eqn:s2b}
t_0^{\alpha} \mid_{\mathcal{H}} =\gwii{\gamma_{1}\gamma^{\alpha}}.
\end{equation}
In particular,  the lifting of $u_i = u^i$ via transformation \eqref{eqn:s2b} and their derivatives with respect to $x:=t_{0}^{1}$ form a coordinate system on the jet space of $\mathcal{H}$. Precise
 definition for functions
$G_{i}^{(2)}$, $G_{ij}^{(2)}$, $P_{ij}^{(2)}$, and $Q_{i}^{(2)}$  can be found in Appendix~\ref{sec:G2}.

To prove that the genus-2 G-function vanishes, it suffices to show that functions $G_{i}^{(2)}$ and
$Q_{i}^{(2)} + \frac{1}{2} P_{ii}^{(2)}$
for all $i$, $G_{ij}^{(2)}$ and
$P_{ij}^{(2)} + P_{ji}^{(2)} $ for $i \neq j$,
are equal to 0. As explained in \cite{LW}, the main idea of the proof is to first represent all these functions in terms of rotation coefficients $r_{ij}$ and
functions $h_{i}$, $u_i$, $\theta_{ij}$, $\Omega_{ij}$ and constant $C$.
Using conditions
$(C1)$, $(C2)$ and $(C3)$, we can then
get rid of $u_i$, $\theta_{ij}$, and $\Omega_{ij}$ and obtain expressions for these functions only
involving $r_{ij}$, $h_{i}$ and $C$. For each of functions $G_{i}^{(2)}$, $G_{ij}^{(2)}$,
$Q_{i}^{(2)} + \frac{1}{2} P_{ii}^{(2)}$,  and
$P_{ij}^{(2)} + P_{ji}^{(2)} $,
the proof of Theorem 0.1 in \cite{LW} shows that
the sum of all terms only containing
$r_{ij}$ and $h_{i}$ is 0. In this paper, we will prove that the sum of all terms containing constant $C$
is also 0.

\subsection{Vanishing of $G_{ij}^{(2)}$}
\label{sec:Gij}

In this subsection, we will prove
 that condition (C2) alone implies that $G_{ij}^{(2)}=0$ for all $i\neq j$.
In fact, by equation \eqref{eq:g=1 1-pt small}, the function $H_i$ in the definition
of the genus-2 G-function (see Appendix~\ref{sec:G2}) has the following form
\begin{align}
H_{i}:= - \frac{1}{2}\sum_{j} v_{ij}r_{ij} =\frac{1}{24}\sum_{k}r_{ik}\frac{h_{i}}{h_{k}}+Ch_{i}^{2}
\label{eq:H}
\end{align}
if condition (C2) is satisfied.
With a slight modification due to the occurrence of constant $C$,
the same computation as in section 2.1 of \cite{LW} shows that
  \begin{align*}
 G_{ij}^{(2)}
 =&-\frac{r_{ij}}{5760h_{i}h_{j}}\big\{{\theta _{ij}}{{{h_j}} \over {{h_i}}} + {\theta _{ji}}{{{h_i}} \over {{h_j}}}- 12r_{ij}^2 - \sum\limits_l {({r_{il}}{r_{jl}}{{{h_i}{h_j}} \over {h_l^2}} - {r_{ij}}{r_{il}}{{{h_j}} \over {{h_l}}} - {r_{ij}}{r_{jl}}{{{h_i}} \over {{h_l}}})}\big\}
-C\frac{r_{ij}^{2}}{120}.
 \end{align*}
The right hand side of this equation is understood as a function obtained by applying
transformation \eqref{eqn:s2b} to
a function on the small phase space.
 Therefore by equation (\ref{eq:g=1 2-pt small}), we have $ G_{ij}^{(2)}=0$.
Most formulas in the proof of Theorem~\ref{thm:main}  will be understood in a similar way.
To make the paper concise, we will not repeat this argument later.

It is not surprising that the vanishing of $G_{ij}^{(2)}$ just follows from a genus-1 condition, i.e.
condition (C2). In fact without assuming condition $(C2)$, one can show that on the small phase space
  \begin{align} \label{eqn:Gijg=1}
  G_{ij}^{(2)}
  ={{{r_{ij}}} \over {240{h_i}{h_j}}}{\phi _{ij}}
  +\frac{r_{ij}^2}{240} \left( \frac{\phi_{i}}{h_{i}^2}
 -\frac{\phi_{j}}{h_{j}^2} \right).
 \end{align}
 So by equations (\ref{eq:g=1 k=1-pt}) and (\ref{eq:g=1 k-pt}), we have $G_{ij}^{(2)}=0$.
 This gives an alternative proof for the vanishing of $G_{ij}^{(2)}$.
 Moreover, equation \eqref{eqn:Gijg=1} and Lemma \ref{lem:C2Idemp} indicate that
 condition (C2) is a rather natural condition for the vanishing of the genus-2 G-function.

\subsection{ Vanishing of $G_{i}^{(2)}$}
\label{sec:Gi}

Since we only need to show that the constant $C$ in condition (C2) does not affect the
final result for the calculation of the genus-2 G-function, in the rest part of the proof
for Theorem~\ref{thm:main}, we will omit those terms whose contributions to the genus-2 G-function
do not contain $C$. To show $G_{i}^{(2)}=0$ for all $i$, condition (C2) is no longer sufficient. We also
need condition (C1). Using similar method as in \cite{LW}, we can compute
 \[ G_{i}^{(2)} = G_{i, 1}^{(2)} + G_{i, 2}^{(2)} \]
and obtain
\begin{align*}
G_{i,1}^{(2)}
=&\sum_{k\neq i}\frac{1}{1920h_{i}^{2}}\big\{{\theta _{ik}}{{{h_k}} \over {{h_i}}}
    + {\theta _{ki}}{{{h_i}} \over {{h_k}}}- 12r_{ik}^2
    - \sum\limits_l {({r_{il}}{r_{kl}}{{{h_i}{h_k}} \over {h_l^2}}
    - {r_{ik}}{r_{il}}{{{h_k}} \over {{h_l}}} - {r_{ik}}{r_{kl}}{{{h_i}} \over {{h_l}}})} \\
& \hspace{30pt}    +48Cr_{ik}h_{i}h_{k}\big\}\frac{\partial_{x}u_{k}}{\partial_{x}u_{i}}
+\sum_{k\neq
i}\frac{\theta_{ik}}{5760h_{i}h_{k}}
+C\frac{r_{ii}}{60}+\{\mbox{some  function  of  $r$  and  $h$}\}
\end{align*}
 and
\begin{align}
G_{i,2}^{(2)}
=&\sum_{k\neq i}\frac{\theta_{ik}}{2880h_{i}h_{k}}-\sum_{k\neq i}\frac{\theta_{ki}h_{k}}{384h_{i}^{3}}
    +\sum_{k\neq i}\frac{7\theta_{ki}}{2880h_{i}h_{k}}
+C\sum_{k}\frac{r_{ik}h_{i}}{120h_{k}}-C\frac{17}{120}r_{ii}
 \nonumber\\&+\{\mbox{some  function  of  $r$  and  $h$}\}
. \label{eqn:Gi2}
\end{align}
In these formulas, the phrase "some  function  of  $r$  and  $h$" means a function which only depends on
$r_{jk}$ and $h_j$ for $j, k =1, \cdots, N$. The contributions of such functions to the genus-2
G-function will not involve the constant $C$.

By equation (\ref{eq:g=1 2-pt small}), the coefficient for
$\frac{\partial_{x}u_{k}}{\partial_{x}u_{i}}$ in $G_{i,1}^{(2)}$ becomes $0$ under condition (C2).
So we have
\begin{align}
G_{i,1}^{(2)}=&\sum_{k\neq
i}\frac{\theta_{ik}}{5760h_{i}h_{k}}
+C\frac{r_{ii}}{60}+\{\mbox{some  function  of  $r$  and  $h$}\}.
\label{eqn:Gi1}
\end{align}

Using conditions (C1) and (C2) , we can then
get rid of first order poles $\theta_{ik}$ and $\theta_{ki}$ in
equations \eqref{eqn:Gi1} and \eqref{eqn:Gi2}.
 A key ingredient in this process is the following lemma:
 \begin{lem} \label{lem:Gi}
Under conditions $(C1)$ and $(C2)$, on the small phase space we have
    \begin{align}
2\sum_{i\neq k}\frac{\theta_{ki}}{h_i}
\,\,\, = \,\,\, &
\sum\limits_i {(7{{h_k} \over {h_i^2}}r_{ik}^2 - 6 \frac{r_{ik}^2}{h_k}
- 2{{h_k^2} \over {h_i^3}}{r_{ii}}{r_{ik}})}
+ \sum\limits_{i,j} {({{h_k^2} \over {{h_i}h_j^2}}{r_{ij}}{r_{jk}}
- {{h_k} \over {{h_i}{h_j}}}{r_{ik}}{r_{jk}})} \nonumber \\
& \hspace{5pt} -24C\sum_{i}\frac{r_{ik}h_{k}^{2}}{h_{i}}
\label{eq:sumtheta}
   \end{align}
for all $k$.
\end{lem}

\noindent{\bf Proof: }
Condition (C1) is equivalent to
\[ \gwii{\ga\gua\gb\gub}|_{\mathcal{H}} = \mbox{constant}. \]
Therefore the derivative of $\gwii{\ga\gua\gb\gub}|_{\mathcal{H}}$ along
$\mathcal {E}_{k}$ must be zero. As explained in \cite{LW}, this results gives the following formula.
\begin{align}
\sum\limits_{i \ne k} ( {{{h_k}} \over {h_i^3}} + {{{h_i}} \over {h_k^3}}){\theta _{ik}} - 2\sum\limits_{i \ne k} {{{{\theta _{ki}}} \over {{h_k}{h_i}}}}
=\{\mbox{some function of $r$ and $h$}\}.
\label{eqn:derO1-O2}
\end{align}
As in \cite{LW}, we can use equation (\ref{eq:g=1 2-pt small}) to
convert the first summation on the left hand side to an expression
similar to the second summation together with some function only
depending on $r_{ij}$, $h_i$ and $C$. Equation \eqref{eq:sumtheta}
is then obtained through straightforward calculations. $\Box$

\begin{rem} \label{rem:Gi}
Note that repeatedly applying equations \eqref{eq:theta},
\eqref{eq:g=1 2-pt small} and $(\ref{eq:sumtheta})$ shows that
\[ \theta_{ij} h_i^p h_j^q \,\, = \,\, \theta_{ij} h_i^{p+2k} h_j^{q-2k} + \{\mbox{some function of $r$, $h$ and $C$}\}. \]
Therefore by equations \eqref{eq:theta} and $(\ref{eq:sumtheta})$,
if $p$ is odd, both $\sum_{i\neq k}h_{i}^{p}\theta _{ik}$
  and $\sum_{i\neq k}h_{i}^{p}\theta _{ki}$
  can be expressed as functions only involving
  $\{ r_{ij} \}$, $\{ h_i \}$ and $C$.
\end{rem}

Using Lemma~\ref{lem:Gi} and Remark~\ref{rem:Gi}, we can get rid of
$\theta_{ik}$ and $\theta_{ki}$ in equations \eqref{eqn:Gi1} and
\eqref{eqn:Gi2} and obtain a formula for $G_{i}^{(2)}=G_{i, 1}^{(2)}
+ G_{i, 2}^{(2)}$ which only involves functions $\{r_{jk} \}$,
$\{h_j \}$ and $C$. In this formula, the part which does not contain
$C$ vanishes as proved in \cite{LW}. So we only need to consider
terms which contain $C$. Note that the constant $C$ appears in the
process of getting rid of $\theta$ terms in $G_{i}^{(2)}$.
By Lemma~\ref{lem:Gi}, the term
$\sum_{k\neq i}\frac{\theta_{ik}}{5760h_{i}h_{k}}$ in $G_{i,1}^{(2)}$ contributes
$-\frac{1}{480}C\sum_{k}\frac{r_{ik}h_{i}}{h_{k}}$. The term $\sum_{k\neq i}\frac{\theta_{ik}}{2880h_{i}h_{k}}$
in $G_{i,2}^{(2)}$ contributes $-\frac{1}{240}C\sum_{k}\frac{r_{ik}h_{i}}{h_{k}}$. By Remark~\ref{rem:Gi},
$-\sum_{k\neq i}\frac{\theta_{ki}h_{k}}{384h_{i}^{3}}$ and $\sum_{k\neq i} \frac{7\theta_{ki}}{2880h_{i}h_{k}}$ in
$G_{i,2}^{(2)}$ contribute $-\frac{1}{32}C\sum_{k}r_{ik}\frac{h_{i}}{h_{k}}+\frac{1}{8}Cr_{ii}$ and
$\frac{7}{240}C\sum_{k}\frac{r_{ik}h_{i}}{h_{k}}$
respectively. Together with terms containing $C$ in equations \eqref{eqn:Gi1} and
\eqref{eqn:Gi2}, the total contribution of all terms containing $C$ in $G_{i}^{(2)}$ is equal to $0$.
This shows that  $G_{i}^{(2)}=0$ under conditions (C1) and (C2).


\subsection{ Vanishing of $P_{ij}^{(2)}+P_{ji}^{(2)}$}
\label{sec:P=0}

In this subsection, we prove $P_{ij}^{(2)}+P_{ji}^{(2)}=0$ for all $i \neq j$ under conditions (C1) and (C2).
After computing all partial derivatives involved in the definition of $P_{ij}^{(2)}$ using
formulas \eqref{r_{ij}} and  \eqref{eq:der of r} and computing functions $H_i$ using equation \eqref{eq:H},
we obtain the following formula for $P_{ij}^{(2)}$:
 \begin{align}
P_{ij}^{(2)}=&\frac{1}{1440} \left( -3 \sum_{k\neq j}\frac{h_{i}h_{j}r_{ik}\theta_{kj}}{h_{k}^{4}}
    -\frac{41r_{ij}\theta_{ij}}{h_{i}^{2}}
+\frac{6r_{ii}h_{j}\theta_{ij}}{h_{i}^{3}}
-3 \sum_{k}\frac{r_{ik}h_{j}\theta_{ij}}{h_{i}^{2}h_{k}} \right.
\nonumber
\\& \left. \hspace{60pt} -4 \sum_{k\neq i}\frac{r_{ij}h_{j}h_{k}\theta_{ik}}{h_{i}^{4}}
- 5 \sum_{k\neq i}\frac{r_{ij}h_{j}\theta_{ki}}{h_{i}^{2}h_{k}} \right)
\nonumber\\&-C\frac{\theta_{ij}h_{j}}{60h_{i}}-C\frac{19r_{ij}^{2}}{30}
-C\frac{r_{ii}r_{ij}h_{j}}{60h_{i}}+\frac{1}{4}C\frac{r_{ij}r_{jj}h_{i}}{h_{j}}
 -\frac{2}{5}C^{2}r_{ij}h_{i}h_{j}
 -\frac{1}{30}C\sum_{k}\frac{r_{ij}r_{jk}h_{i}}{h_{k}}
\nonumber\\&-C\sum_{k}\frac{r_{ik}r_{jk}h_{i}h_{j}}{60h_{k}^{2}}
+\{\mbox{some function of $r$ and $h$}\}.
\label{eqn:Pijp}
\end{align}
As explained before we only need to consider terms whose contributions to the final calculation of
$P_{ij}^{(2)}$ may contain the constant $C$. The omitted terms in the above formula will not contribute
constant $C$ in the final calculation.

When getting rid of terms containing $\theta$  in equation \eqref{eqn:Pijp}, new terms containing constant $C$
will emerge.
We can get rid of the fifth and sixth terms on the right hand side
 of equation \eqref{eqn:Pijp} directly using Lemma~\ref{lem:Gi} and Remark~\ref{rem:Gi}.
They will contribute
 \[ -\frac{2}{15}C\frac{r_{ij}r_{ii}h_{j}}{h_{i}}-\frac{1}{120}Cr_{ij}h_{j}
\sum_{k}\frac{r_{ik}}{h_{k}}. \]
To deal with other terms containing $\theta$ in equation \eqref{eqn:Pijp},
 we need consider the sum
 \[ P_{ij}^{(2)}+P_{ji}^{(2)}. \]
We start with the calculation for the first term on the right hand side of equation \eqref{eqn:Pijp}.
 By equation (\ref{eq:g=1 2-pt small}), we have
         \begin{align}
\sum_{k\neq j}\frac{h_{i}h_{j}r_{ik}\theta_{kj}}{h_{k}^{4}}
= - \frac{h_{i}}{h_{j}}\sum_{k\neq j}\frac{r_{ik}\theta_{jk}}{h_{k}^{2}}
-48C\sum_{k\neq j}\frac{r_{ik}r_{jk}h_{i}h_{j}}{h_{k}^{2}}+\{\mbox{some  function  of  $r$  and  $h$}\}. \label{P}
\end{align}
Using the similar method as in \cite{LW}, we can show that condition $(C2)$ implies the following formula
on the small phase space
\begin{align}
& {\theta _{ik}}{r_{jk}}\frac{{{h_j}}}{{{h_i}}} + {\theta _{ki}}{r_{ij}}\frac{{{h_j}}}{{{h_k}}}
 + {\theta _{jk}}{r_{ik}}\frac{{{h_i}}}{{{h_j}}} + {\theta _{kj}}{r_{ij}}\frac{{{h_i}}}{{{h_k}}}
 + {\theta _{ij}}r_{jk}\frac{{{h_k}{}}}{{{h_i}}} + {\theta _{ji}}{r_{ik}}\frac{{{h_k}}}{{{h_j}}}
 \nonumber\\&+72Cr_{ik}r_{jk}h_{i}h_{j}+72Cr_{ij}r_{ik}h_{j}h_{k}
 +72Cr_{ij}r_{jk}h_{i}h_{k}
                            \nonumber
 \\=&\{\mbox{some  function  of  $r$  and  $h$}\}
 \label{g=1 3-pt jian}
\end{align}
  for  $i\neq j\neq k$.
Multiplying both sides of equation (\ref{g=1 3-pt jian}) by $\frac{1}{h_{k}^{2}}$
and taking sum over $k$ for $k\notin \{i,j\}$,  we  obtain
   \begin{align}
   & \frac{h_{i}}{h_{j}}\sum_{k\neq j}\frac{r_{ik}\theta_{jk}}{h_{k}^{2}}+ \frac{h_{j}}{h_{i}}\sum_{k\neq j}\frac{r_{jk}\theta_{ik}}{h_{k}^{2}}
  \nonumber\\=&  - \sum_{k}{\theta _{ij}}r_{jk}\frac{{{1}}}{{{h_i}{h_k}}}
 - \sum_{k}{\theta _{ji}}{r_{ik}}\frac{{1}}{{{h_j}{h_k}}}
+2{\theta _{ij}}{r_{jj}}\frac{{1}}{{{h_i}{h_j}}}
  +2{\theta _{ji}}{r_{ii}}\frac{{1}}{{{h_i}{h_j}}}
-36Cr_{ij}h_{j}\sum_{k}\frac{r_{ik}}{h_{k}}
 \nonumber\\&- 36 Cr_{ij}h_{i}\sum_{k}\frac{r_{jk}}{h_{k}}
-72C\sum_{k}\frac{1}{h_{k}^{2}}{r_{ik}}{r_{jk}}h_{i}h_{j}
+96 Cr_{ii}{r_{ij}}\frac{h_{j}}{h_{i}}
+96 Cr_{jj}{r_{ij}}\frac{h_{i}}{h_{j}}
+48 Cr_{ij}^{2}
   \nonumber\\&+\mbox{\{some function of $r$ and $h$\}}. \label{P1}
   \end{align}
In deriving this formula, we have used equation \eqref{eq:g=1 2-pt small} and the fact that both
$\sum_{k \neq i} \frac{\theta_{ki}}{h_k^3}$ and
$\sum_{k \neq j} \frac{\theta_{kj}}{h_k^3}$
can be expressed as functions of $\{r_{kl} \}$, $\{h_k\}$ and $C$.

First applying equation \eqref{P} to the first term of $P_{ij}^{(2)}$ on the right hand side of
equation \eqref{eqn:Pijp}, then adding to the corresponding term in $P_{ji}^{(2)}$, we can apply
equation \eqref{P1} to replace these terms by other terms which
still contain $\theta$. After this operation, all the terms in $P_{ij}^{(2)} +P_{ji}^{(2)}$
containing $\theta$
can be grouped into several pairs for which equation \eqref{eq:g=1 2-pt small} can be applied.
Using equation \eqref{eq:g=1 2-pt small}, we can then obtain a formula for
 $P_{ij}^{(2)} +P_{ji}^{(2)}$
which only depends on $\{r_{kl}\}$, $\{h_{k}\}$ and $C$. For example, if we just get rid of
terms in $P_{ij}^{(2)} +P_{ji}^{(2)}$ which contain $\theta$ but do not contain $C$
using the above procedure, we obtain
the following formula
  \begin{align*}
&P_{ij}^{(2)}+P_{ji}^{(2)}
\\=&
 - C\frac{{\theta _{ij} h_j }}{{60h_i }} - C\frac{{\theta _{ji} h_i }}{{60h_j }} - \frac{1}{{60}}Cr_{ij} h_i \sum\limits_k {r_{jk} \frac{1}{{h_k }}}  - \frac{1}{{60}}Cr_{ij} h_j \sum\limits_k {r_{ik} \frac{1}{{h_k }}}  + \frac{1}{{60}}C\sum\limits_k {\frac{{r_{ik} r_{jk} h_i h_j }}{{h_k^2 }}}
\\&
 + C\frac{1}{5}r_{ij}^2  - \frac{4}{5}C^2 r_{ij} h_i h_j
 +\mbox{\{some function of $r$ and $h$\}}
.
\end{align*}
Applying equation \eqref{eq:g=1 2-pt small} again to the first two terms on the right hand side,
we see that the sum of all terms containing $C$ in this expression is 0.
Together with the result in \cite{LW}, this proves that
$P_{ij}^{(2)}+P_{ji}^{(2)}=0$ for $i \neq j$.

\subsection{ Vanishing of $\frac{1}{2}P_{ii}^{(2)}+Q_{i}^{(2)}$}
\label{sec:P+Q=0}

In this subsection, we prove $\frac{1}{2}P_{ii}^{(2)}+Q_{i}^{(2)}=0$ for all $i$.
For this purpose, we also need condition (C3) in addition to conditions (C1) and (C2).
Using the same method as in \cite{LW}, we can express
$\frac{1}{2}P_{ii}^{(2)}+Q_{i}^{(2)}$  as a function of
$\{ \Omega_{ij} \}$,  $\{\theta_{jk} \}$, $\{ v_{jk} \}$,  $\{ r_{jk} \}$, $\{ h_{j} \}$ and $C$.
After getting rid of part of $\{\theta_{jk} \}$ terms using equations \eqref{eq:theta}, \eqref{eq:g=1 2-pt small},
Lemma~\ref{lem:Gi} and Remark~\ref{rem:Gi}, we have
           \begin{align}
&5760\{\frac{1}{2}P_{ii}^{(2)}+Q_{i}^{(2)}\}\nonumber
\\=&- 5\sum_{k\neq i} {\frac{{{\widetilde{\Omega} _{ik}}}}{{h_i^2}}}  - 100\sum\limits_{k \ne i} {\frac{{{r_{ik}}{\theta _{ik}}}}{{h_i^2}}}  + 4\sum\limits_{k \ne i} {\frac{{{r_{kk}}{h_i}{\theta _{ik}}}}{{h_k^3}}}  - 6\sum\limits_{k \ne i} {\frac{{h_i^2{r_{ik}}{\theta _{ki}}}}{{h_k^4}}}
  + 5\sum_{j}\sum_{k\neq i}{\frac{{{r_{jk}}{\theta _{ik}}}}{{{h_i}{h_j}}}}
  \nonumber\\&- 2\sum_{k}\sum_{j\neq i}{\frac{{{r_{jk}}{h_i}{\theta _{ij}}}}{{{h_k}h_j^2}}}  + 40\sum_{k}\sum_{j\neq i}{\frac{{{r_{ik}}{v_{jk}}{\theta _{ij}}}}{{h_i^2}}}
   -8064{C^2}{r_{ii}}h_i^2 - 48Cr_{ii}^2 +4608{C^3}h_i^4 - 856C\sum\limits_k {r_{ik}^2}
   \nonumber\\&- 48C\sum\limits_k {\frac{{r_{ik}^2}}{{h_k^2}}} h_i^2 + 456C{r_{ii}}{h_i}\sum\limits_k {\frac{{{r_{ik}}}}{{{h_k}}}}
   + 576{C^2}\sum\limits_k {\frac{{{r_{ik}}}}{{{h_k}}}} h_i^3 + 168C\sum_{j,k} {\frac{{{r_{ik}}{r_{jk}}}}{{{h_j}}}} {h_i}
   \nonumber\\&- 108C\sum_{j,k} {\frac{{{r_{ik}}{r_{ij}}}}{{{h_k}{h_j}}}} h_i^2
  + 168C\sum\limits_{k \ne i} {\frac{{{\theta _{ki}}{h_i}}}{{{h_k}}}}  + 48C\sum\limits_{k \ne i} {\frac{{{\theta _{ik}}{h_i}}}{{{h_k}}}}  + 456C\sum\limits_{k \ne i} {\frac{{{\theta _{ik}}{h_k}}}{{{h_i}}}}
  \nonumber\\& +\{\mbox{some function of $r$ and $h$}\},
   \label{Q+P jian}
           \end{align}
where $\widetilde{\Omega}_{ij}=\Omega_{ij}\{\frac{h_{i}}{h_{j}}-
\frac{h_{j}}{h_{i}}\}$.
Condition (C3) is needed to deal with the term
\[ 40\sum_{k}\sum_{j\neq i} {\frac{{{r_{ik}}{v_{jk}}{\theta _{ij}}}}{{h_i^2}}} \]
in this expression.
 Here we also need to use the formula of $F_2$ obtained in \cite{L07}.
 Following the approach described in the proof of Lemma 2.5 in \cite{LW},  we obtain
\begin{align}
&80\sum_{j\neq i} {\sum_{k} {{\theta _{ij}}{r_{ik}}{v_{jk}}\frac{1}{{h_i^2}}} }
 \nonumber
\\=& - 15\frac{1}{{h_i^2}}\sum_{j\neq i} {{{\widetilde \Omega }_{ji}}}  - 5\sum_{j \neq i}{{{\widetilde \Omega }_{ji}}\frac{1}{{h_j^2}}}  - 24\sum_{j\neq i} {{\theta _{ji}}{r_{jj}}\frac{{{h_i}}}{{h_j^3}}}
 - 400\sum_{j\neq i} {{\theta _{ji}}{r_{ji}}\frac{1}{{h_i^2}}} + 22\sum_{j\neq i} {\sum\limits_k {{\theta _{ji}}{r_{jk}}\frac{1}{{{h_k}{h_i}}}} }
 \nonumber\\& - 3456Cr_{ii}^2 + 13824{C^2}{r_{ii}}h_i^2 - 9216{C^3}h_i^4
  - 864C{r_{ii}}{h_i}\sum\limits_j {\frac{{{r_{ij}}}}{{{h_j}}}}  + 6128C\sum\limits_k {r_{ik}^2}
   \nonumber\\&- 1152{C^2}\sum\limits_j {{r_{ij}}\frac{{h_i^3}}{{{h_j}}}}
 + 264C\sum\limits_k {\sum\limits_j {{r_{ik}}{r_{ij}}\frac{{{h_i}}}{{{h_k}}}\frac{{{h_i}}}{{{h_j}}}} }  - 528C\sum\limits_j {\sum\limits_k {{r_{jk}}} {r_{ij}}\frac{{{h_i}}}{{{h_k}}}} \nonumber\\& - 1104C\sum_{j\neq i} {\{ {\theta _{ij}}\frac{{{h_j}}}{{{h_i}}} + {\theta _{ji}}\frac{{{h_i}}}{{{h_j}}}\} }  + 1152C\sum_{j\neq i}{{\theta _{ji}}\frac{{{h_i}}}{{{h_j}}}}
 +\{\mbox{some function of $r$ and $h$}\}
 \label{TE der g=2jian}
 \end{align}
for all $i$ when conditions (C1), (C2) and (C3) are satisfied.
Moreover under conditions (C1) and (C2), the formulas in Lemma 2.6 and Lemma 2.7 in \cite{LW} can be modified in the following way
\begin{align}
0=&{\widetilde\Omega}_{ij} + {\theta _{ij}}\{ 20{r_{ij}} + 4{r_{ii}}\frac{{{h_j}}}{{{h_i}}}
- \sum\limits_k {} ({r_{ik}}\frac{{{h_j}}}{{{h_k}}} + {r_{jk}}\frac{{{h_i}}}{{{h_k}}})\}
+ 2\sum\limits_{k \ne i} {{\theta _{ik}}} {r_{jk}}\frac{{{h_j}}}{{{h_i}}}
\nonumber\\& - 144Cr_{ij}^2h_i^2 - 96C{r_{ij}}{r_{ii}}{h_i}{h_j} + 24C{r_{ij}}{h_i}{h_j}\sum\limits_k {\frac{{{r_{ik}}}}{{{h_k}}}} {h_i} + 48C\sum\limits_k {{r_{jk}}} {r_{ik}}{h_i}{h_j}
 \nonumber\\&- 48C{\theta _{ij}}{h_i}{h_j}
   +\{\mbox{some function of $r$ and $h$}\}
   \label{g=13-pt2jian}
\end{align}
for all $i\neq j$, and
     \begin{align}
   &2\sum_{k\neq i}{\theta _{ik}}{r_{ik}}\frac{1}{{h_i^2}}
   - \sum_{j}\sum_{k\neq i}{\theta _{ik}}{r_{jk}}\frac{{1}}{{{h_i}{h_j}}} + 2\sum_{k\neq i} {\theta _{ik}}{r_{kk}}\frac{{{h_i}}}{{h_k^3}}
   \nonumber
   \\=&240Cr_{ii}^2 - 144C\sum\limits_k {r_{ik}^2}  - 72C\sum\limits_k {r_{ik}^2\frac{{h_i^2}}{{h_k^2}}}
    - 96C{r_{ii}}{h_i}\sum\limits_k {\frac{{{r_{ik}}}}{{{h_k}}}}  + 48C{h_i}\sum\limits_{j,k} {\frac{{{r_{ij}}{r_{jk}}}}{{{h_k}}}}
      \nonumber\\&+ 12Ch_i^2\sum\limits_k {\sum\limits_j {\frac{{{r_{ij}}{r_{ik}}}}{{{h_j}{h_k}}}} }
  +\{\mbox{some function of $r$ and $h$}\} \label{g=1 3-pt3}
     \end{align}
 for all $i$. These formulas should be understood as formulas on the small phase space.

After plugging equation \eqref{TE der g=2jian} into equation \eqref{Q+P jian} and using
equation \eqref{g=13-pt2jian} to get rid of ${\widetilde\Omega}_{ij}$,
 we get
\begin{align}
&5760\{\frac{1}{2}P_{ii}^{(2)}+Q_{i}^{(2)}\}
\nonumber\\=&16\frac{1}{{h_i^2}}\sum\limits_{k \ne i} {{r_{ik}}{\theta _{ik}}}  - 8\sum\limits_{\mathop {k,l}\limits_{l \ne i} } {{r_{kl}}{\theta _{il}}\frac{1}{{{h_k}}}\frac{1}{{{h_i}}}}  +16\sum\limits_{k \ne i} {\frac{{{r_{kk}}{h_i}{\theta _{ik}}}}{{h_k^3}}}
 + 1152C\sum\limits_k {r_{ik}^2} - 1920Cr_{ii}^2
 \nonumber\\& + 576Ch_i^2\sum\limits_j {r_{ij}^2\frac{1}{{h_j^2}}}  + 768C{r_{ii}}{h_i}\sum\limits_k {\frac{{{r_{ik}}}}{{{h_k}}}}
- 384C{h_i}\sum\limits_j {\sum\limits_k {{r_{ij}}{r_{jk}}\frac{1}{{{h_k}}}} }
\nonumber\\&   - 96C{h_i}{h_i}\sum\limits_k {\sum\limits_j {{r_{ik}}{r_{ij}}\frac{1}{{{h_j}{h_k}}}} }
 +\mbox{\{some function of $r$ and $h$\}}.
\label{eqn:PQ3theta}
\end{align}
In this process, we have used equation \eqref{eq:g=1 2-pt small}, Lemma~\ref{lem:Gi} and
Remark~\ref{rem:Gi} to get rid of terms  containing $\theta$ as many as possible.
We also notice that terms containing $C^2$ and $C^3$ are all canceled in this process.
For example, the term $9216C^{3}h_{i}^{4}$ in equation \eqref{Q+P jian}
is canceled when plugging
equation \eqref{TE der g=2jian} into equation \eqref{Q+P jian}.

Plugging equation \eqref{g=1 3-pt3} into equation \eqref{eqn:PQ3theta}, we see that
$\frac{1}{2}P_{ii}^{(2)} + Q_i^{(2)}$ does not depend on the constant $C$. Together with the
proof in \cite{LW}, this shows that
  $\frac{1}{2}P_{ii}^{(2)} + Q_i^{(2)}=0$ for all $i$.
Theorem~\ref{thm:main} is thus proved.
$\Box$

\section{Proof of Theorem~\ref{thm:ADEG} }
\label{sec:pfcor}
In this section we will prove that conditions $(C1)-(C3)$ are satisfied for $\mathbb{P}^{1}$-orbifolds of $ADE$ type.
We can then use Theorem~\ref{thm:main} to prove Theorem~\ref{thm:ADEG}.

Orbifold Gromov-Witten theory was introduced in
\cite{CR} in symplectic geometric setting. Algebraic geometric
treatments for such invariants were given in \cite{AGV}. In this paper we only need some
basic properties of the Gromov-Witten theory for
$\mathbb{P}^{1}$-orbifolds which can also be found in \cite{KS}. Let
$\mathcal {X}=\mathbb{P}^{1}_{o_{1},o_{2},o_{3}}$ be the orbifold
$\mathbb{P}^{1}$ with three orbifold points such that the $i$-th
orbifold point has isotropy group $\mathbb{Z}/o_{i}\mathbb{Z}$,
where $o_i$ are positive integers. The Chen-Ruan cohomology group of
$\mathcal {X}$ has the following form:
\begin{align*}
H^{*}_{CR}(\mathcal {X})
=\mathbb{C}[\Delta_{01}]\oplus\Big(\oplus_{i=1}^{3}\oplus_{j=1}^{o_{i}-1}\mathbb{C}[\Delta_{ij}]\Big)
\oplus\mathbb{C}[\Delta_{02}]
\end{align*}
where $\Delta_{01}$ is the identity and $\Delta_{02}=\omega$ is the hyperplane class of the underlying $\mathbb{P}^{1}$.
The classes
$\Delta_{ij}$ with $1\leq i\leq3,1\leq j\leq o_{i}-1$ are in one-to-one correspondence with the twisted sectors and
we define $\Delta_{ij}$ to be the unit in the cohomology of the corresponding twisted sector.
The complex degrees of these classes are
\begin{align*}
\mbox{deg} \, \Delta_{01}=0,\quad\mbox{deg} \, \Delta_{02}=1,
\end{align*}
\begin{align*}
\mbox{deg} \, \Delta_{ij}=\frac{j}{o_{i}} \,\,\,\,\,\, \mbox{for}\quad 1\leq i\leq3, \quad 1\leq j\leq o_{i}-1.
\end{align*}
As pointed in \cite{KS}, the orbifold Poincare pairing takes the form
\begin{align*}
<\Delta_{i_{1},j_{1}},\Delta_{i_{2},j_{2}}>
=\begin{cases}
(\delta_{i_{1}, i_{2}}\delta_{j_{1}+j_{2},o_{i_{1}}})/o_{i_{1}}&\quad\quad \text{if $i_{1}+i_{2}\neq 0$},\\
\delta_{j_{1}+j_{2},3}& \quad\quad\text{if $i_{1}=i_{2}=0.$}
\end{cases}
\end{align*}
The genus-0 3-point degree-0 Gromov-Witten invariants are given by the following formula
\begin{align*}
& <\Delta_{i_{1},j_{1}},\Delta_{i_{2},j_{2}},\Delta_{i_{3},j_{3}}>_{0,0}  \\
= & \begin{cases}
1/o_{i_{1}}&\quad\quad \text{if $i_{1}=i_{2}=i_{3}\in \{1,2,3\}, j_{1}+j_{2}+j_{3}=o_{i_{1}}$},\\
<\Delta_{i_{2}j_{2}},\Delta_{i_{3}j_{3}}>& \quad\quad\text{if $(i_{1},j_{1})=(0,1)$},\\
0& \quad\quad\text{otherwise}.
\end{cases}
\end{align*}

Let $\overline{\mathcal {M}}_{g,n}(\mathcal {X},\beta)$ be the moduli space of stable maps from genus-$g$, $n$-pointed
orbi-curves to $\mathcal {X}$ with the push forward of the fundamental class equal to $\beta$, which lies in the
Mori cone of the homological classes of effective 1-cycles. In our case, since $\mathcal {X}$ is 1-dimensional, we can
write $\beta=d[\mathcal {X}]$ for $d\geq 0$, where $[\mathcal {X}]$ denotes the fundamental class of orbifold
$\mathcal {X}$ . The virtual dimension of $\overline{\mathcal {M}}_{g,n}(\mathcal {X},\beta)$ is
\begin{align*}
\mbox{vir dim}_{\mathbb{C}}\overline{\mathcal {M}}_{g,n}(\mathcal {X},\beta)
&=(3-1)(g-1)+n+<c_{1}(T\mathcal {X}),\beta> \\
&=2g-2+n+d(\sum_{i=1}^{3}\frac{1}{o_{i}}-1)
\end{align*}
The descendant Gromov-Witten invariants of $\mathcal{X}$ satisfy the divisor equation (cf. \cite{AGV}):
\begin{align}
\gwigd{\tau_{n_{1}}(\gamma_{\alpha_{1}})\cdot\cdot\cdot
\tau_{n_{k}}(\gamma_{\alpha_{k}})\omega}
= & \,\, d \, \gwigd{\tau_{n_{1}}(\gamma_{\alpha_{1}})\cdot\cdot\cdot
\tau_{n_{k}}(\gamma_{\alpha_{k}})}  \nonumber \\
& + \sum_{i=1}^k \gwigd{\tau_{n_{1}}(\gamma_{\alpha_{1}})\cdots \tau_{n_{i}-1}(\omega \cup \gamma_{\alpha_{1}}) \cdots
\tau_{n_{k}}(\gamma_{\alpha_{k}})}
\label{divisor eq}
\end{align}
for all $k\geq1$, $d\geq0$, $n_{i}\geq0$, for $1\leq i\leq k$ and $\gamma_{\alpha_{1}},...,\gamma_{\alpha_{k}}$
in Chen-Ruan cohomology group of $\mathbb{P}^{1}_{o_{1},o_{2},o_{3}}$.
By definition, a $\mathbb{P}^{1}$-orbifold of $ADE$ type is $\mathbb{P}^{1}_{o_{1},o_{2},o_{3}}$ with property
\[\sum_{i=1}^{3}\frac{1}{o_{i}}-1>0.\]
\begin{thm} \label{P1inv}
For $\mathbb{P}^{1}$-orbifolds of $ADE$ type, we have
\begin{align*}
 (i)& \quad \gwid{ {\gamma _\alpha }{\gamma ^\alpha }{\gamma _\beta }{\gamma ^\beta }\gamma_{\alpha_{1}}
 \cdot\cdot\cdot \gamma_{\alpha_{k}}} =0  \hspace{10pt}\mbox{for any fixed}
 \hspace{5pt}\alpha \hspace{5pt}\mbox{and}\hspace{5pt} \beta,
\\
 (ii)&\quad   \gwioned{ \gamma_{\alpha_{1}}\cdot\cdot\cdot \gamma_{\alpha_{k}}} =\begin{cases}
-\frac{1}{24}& \text{if \hspace{5pt}$k=1,d=0,\gamma_{\alpha_{1}}=\omega$,}\\
0& \text{\hspace{5pt}otherwise,}
\end{cases} \\
 (iii)&\quad \gwitwod{ \gamma_{\alpha_{1}}\cdot\cdot\cdot \gamma_{\alpha_{k}}} =0 \quad \mbox{and} \quad
\gwitwod{ \tau_{1}(\gamma_{\alpha_{1}})\gamma_{\alpha_{2}}\cdot\cdot\cdot \gamma_{\alpha_{k}}} =0,
 \end{align*} for all $k\geq1$, $d\geq0$ and any $\gamma_{\alpha_{1}},...,\gamma_{\alpha_{k}}$ in Chen-Ruan cohomology group of $\mathbb{P}^{1}_{o_{1},o_{2},o_{3}}$.
\end{thm}
{\bf Proof}: Consider the correlator $\gwid{ {\gamma _\alpha }{\gamma ^\alpha }{\gamma _\beta }{\gamma ^\beta }\gamma_{\alpha_{1}}
 \cdot\cdot\cdot \gamma_{\alpha_{k}}}$.
The orbifold poincare paring on Chen-Ruan cohomology $H^{*}_{CR}(\mathcal {X})$
implies that the sum of the complex degrees of $\gamma_{\alpha}$ and $\gamma^{\alpha}$ (also the sum of the
complex degrees of $\gamma_{\beta}$ and $\gamma^{\beta}$) must be equal to $1$, i.e. the complex dimension of
$\mathcal {X}=\mathbb{P}^{1}_{o_{1},o_{2},o_{3}}$. Since the complex degree of all cohomology classes
are  $\leq 1$, the sum of degrees of all cohomology classes in this correlator is
$\leq 2+k$. On the other hand, the complex dimension of the virtual fundamental class equals to
$k+2+d(\sum_{i=1}^{3}\frac{1}{o_{i}}-1)$. Thus for $d\geq 1$, this correlator must be 0 for the dimension reason.
For $d=0$, the dimension constraint requires all $\gamma_{\alpha_{1}},..., \gamma_{\alpha_{k}}$ equal to
$\omega$. The divisor equation \eqref{divisor eq} implies that this correlator is $0$ for $d=0$.
This proves $(i)$.
Similarly, by dimension constraint, divisor equation and the fact $\gwioneo{\omega}=-\frac{1}{24}$ for
$\mathbb{P}^{1}_{o_{1},o_{2},o_{3}}$, we obtain $(ii)$. Part $(iii)$ just follows from dimension constraint.
$\Box$

Part $(i)$ and part $(ii)$ of  Theorem~\ref{P1inv} implies the following
\begin{cor} \label{cor:P1ADEg=1G}
Let $\mathcal {X}$ be a $\mathbb{P}^{1}$ orbifold of $ADE$ type.
Genus-0 and genus-1 primary Gromov-Witten invariants of $\mathcal {X}$
satisfy the following identities:
\begin{align} \label{eqn:C1P1}
\sum_{\alpha,\beta}\gwii{\gamma_{\alpha}\gamma^{\alpha}
\gamma_{\beta}\gamma^{\beta}}|_{\mathcal{H}}= \mbox{constant}
\end{align}
and
\[F_1 \mid_{\mathcal{H}}=-\frac{1}{24}t^{N} + \mbox{constant}, \]
where $F_1 \mid_{\mathcal{H}}$ is the generating function for genus-1 primary invariants and
$t^{N}$ is the flat coordinate on the small phase space $\mathcal{H}= H^{*}_{CR}(\mathcal {X})$
corresponding to the cohomology class $\omega$.
\end{cor}
Note that by Lemma 3.1 of \cite{LW}, after the transformation \eqref{eqn:s2b},
the function
\[ \sum_{\alpha,\beta}\gwii{\gamma_{\alpha}\gamma^{\alpha}
\gamma_{\beta}\gamma^{\beta}}|_{\mathcal{H}} \]
is equal to $O_1 -O_2$ where
\[O_{1}=\gwii{\gamma_{\alpha}\gamma_{\alpha'}\gamma_{\beta}
        \gamma_{\beta'}}(M^{-1})^{\alpha\alpha'} (M^{-1})^{\beta\beta'}, \]
\[ O_{2}=\gwii{\gamma_{\alpha}\gamma_{\beta}\gamma_{\rho}}
    \gwii{\gamma_{1}\gamma_{\alpha'}\gamma_{\beta'}\gamma_{\rho'}}
    (M^{-1})^{\alpha\alpha'}(M^{-1})^{\beta\beta'}(M^{-1})^{\rho\rho'} \]
and entries of the matrix $M$ are defined by
\[
 M_{\mu\rho}=\gwii{\gamma_{1}\gamma_{\mu}\gamma_{\rho}}
\]
 for any $\mu$ and $\rho$.
Therefore the special case of Corollary~\ref{cor:P1ADEg=1G} for
$\mathbb{P}^{1}$ orbifolds of $AD$ type
 was also proved in \cite{DLZ} via case by
case studies. For type-$E$ $\mathbb{P}^{1}$ orbifolds, the result of
this corollary was conjectured in \cite{DLZ} with a precise constant
in equation \eqref{eqn:C1P1}. Corollary \ref{cor:P1ADEg=1G} in
particular solves conjecture 3.20 in \cite{DLZ}
 up to a constant. This constant is not important for our proof of
 Theorem \ref{thm:ADEG}. Moreover our proof of Corollary \ref{cor:P1ADEg=1G} is a unified approach for all
 $\mathbb{P}^{1}$ orbifolds of $ADE$ type which is much simpler than the arguments in
 \cite{DLZ} for $\mathbb{P}^{1}$ orbifolds of $AD$ type.

Theorem~\ref{P1inv} implies that conditions (C1)--(C3) are satisfied for $\mathbb{P}^{1}$-orbifolds of $ADE$ type.
Therefore Theorem~\ref{thm:ADEG} follows from Theorem~\ref{thm:main}.

\vspace{40pt}

\appendix
\centerline {\bf \Large Appendix}
\section{The genus-2 G-function}
\label{sec:G2}

In this appendix, we give the precise definition of the genus-2 G-function $G^{(2)}$ following \cite{DLZ}.
Write
\begin{align}
                   G^{(2)}
                   =&\sum_{i}G_{i}^{(2)}(u,u_{x})u_{xx}^{i}
                   +\sum_{i\neq j}G_{ij}^{(2)}(u)\frac{(u_{x}^{j})^{3}}{u_{x}^{i}}
                   \nonumber\\&+\frac{1}{2}\sum_{i,j}P_{ij}^{(2)}(u)u_{x}^{i}u_{x}^{j}
                   +\sum_{i}Q_{i}^{(2)}(u)(u_{x}^{i})^{2}.
                    \label{eqn:g2G}
\end{align}
Let $\gamma_{ij}$ be the rotation coefficient on the small phase space as defined in \cite{D}. Note that
$\gamma_{ii}=0$ which is different from our definition of $r_{ii}$ in Section~\ref{sec:prel}. For $i \neq j$,
$\gamma_{ij}$ is equal to our definition of $r_{ij}$  restricted to the small phase space. Define
\begin{align*}
H_{i}:=\frac{1}{2}\sum_{j\neq i}u_{ij}\gamma_{ij}^{2}
\end{align*}
where $u_{ij}:=u_{i}-u_{j}$.
Then the function $ G_{i}^{(2)}$ can be defined as
\begin{align*}
   G_{i}^{(2)}=G_{i,1}^{(2)}+G_{i,2}^{(2)}
  \end{align*}
with
\begin{align*}
G_{i,1}^{(2)}=&\frac{\partial_{x}h_{i}H_{i}}{60u_{i,x}h_{i}^{3}}
-\frac{7\partial_{i}h_{i}\partial_{x}h_{i}}{5760u_{i,x}h_{i}^{4}}
+\sum_{k}\Bigg(\frac{\gamma_{ik}H_{k}}{120h_{i}h_{k}}\frac{u_{k,x}}{u_{i,x}}
-\frac{\gamma_{ik}\partial_{x}h_{i}}{5760h_{i}^{2}h_{k}u_{i,x}}
-\frac{\gamma_{ik}\partial_{k}h_{k}u_{k,x}}{1152h_{i}h_{k}^{2}u_{i,x}}
\\&+\frac{\partial_{i}\gamma_{ik}h_{k}u_{k,x}}{1920u_{i,x}h_{i}^{3}}
+\frac{\partial_{x}\gamma_{ik}}{5760u_{i,x}h_{i}h_{k}}
+\frac{\partial_{k}\gamma_{ik}u_{k,x}}{2880h_{i}h_{k}u_{i,x}}
-\frac{7\gamma_{ik}^{2}u_{k,x}}{1152h_{i}^{2}u_{i,x}}\Bigg)
-\sum_{k,l}\frac{u_{k,x}h_{k}\gamma_{il}\gamma_{kl}}{1920u_{i,x}h_{i}h_{l}^{2}}
\end{align*}
and
\begin{align*}
G_{i,2}^{(2)}=&-\frac{3\partial_{i}h_{i}H_{i}}{40h_{i}^{3}}+\frac{19(\partial_{i}h_{i})^{2}}{2880h_{i}^{4}}
+\sum_{k}\Bigg(\frac{\gamma_{ik}H_{i}}{120h_{i}h_{k}}+\frac{7\gamma_{ik}H_{k}}{120h_{i}h_{k}}-\frac{4\gamma_{ik}\partial_{i}h_{i}}{5760h_{i}^{2}h_{k}}
-\frac{7\gamma_{ik}\partial_{k}h_{k}}{2880h_{i}h_{k}^{2}}\\&+\frac{\gamma_{ik}\partial_{k}h_{k}}{384h_{i}^{3}}-\frac{\partial_{k}\gamma_{ik}h_{k}}{384h_{i}^{3}}
+\frac{\partial_{i}\gamma_{ik}}{2880h_{i}h_{k}}+\frac{7\partial_{k}\gamma_{ik}}{2880h_{i}h_{k}}+\frac{\gamma_{ik}h_{i}\partial_{k}h_{k}}{2880h_{k}^{4}}
-\frac{19\gamma_{ik}^{2}}{720h_{i}^{2}}+\frac{\gamma_{ik}^{2}}{1440h_{k}^{2}}\Bigg)
\\&-\sum_{k,l}\frac{h_{i}\gamma_{il}\gamma_{kl}}{2880h_{k}h_{l}^{2}}.
\end{align*}
Other functions in equation \eqref{eqn:g2G} are defined in the following way:
   \begin{align*}
G_{ij}^{(2)} =&  - \frac{{\gamma _{ij}^2{H_j}}}{{120h_j^2}} + \frac{{\gamma _{ij}^3}}{{480{h_i}{h_j}}} - \frac{{{\gamma _{ij}}}}{{5760}}(\frac{{{\partial _i}{\gamma _{ij}}}}{{h_i^2}} + \frac{{{\partial _j}{\gamma _{ij}}}}{{h_j^2}}) + \frac{{\gamma _{ij}^2}}{{5760}}(\frac{{{\partial _i}{h_i}}}{{h_i^3}} + \frac{{3{\partial _j}{h_j}}}{{h_j^3}})
\\& + \sum\limits_k {} (\frac{{{\gamma _{ij}}{\gamma _{ik}}{\gamma _{jk}}}}{{5760h_k^2}} + \frac{{\gamma _{ij}^2}}{{5760{h_k}}}(\frac{{{\gamma _{jk}}}}{{{h_j}}} - \frac{{{\gamma _{ik}}}}{{{h_i}}})),
  \end{align*}
   \begin{align*}
P_{ij}^{(2)} =&  - \frac{{2{\gamma _{ij}}{H_i}{H_j}}}{{5{h_i}{h_j}}} + \frac{{{\gamma _{ij}}{\partial _j}{h_j}{H_i}}}{{20{h_i}h_j^2}} + \frac{{{\gamma _{ij}}{h_i}{\partial _j}{h_j}{H_j}}}{{20h_j^4}} - \frac{{19\gamma _{ij}^2{H_j}}}{{30h_j^2}} - \frac{{{\partial _i}{\gamma _{ij}}{H_j}}}{{60{h_i}{h_j}}}
\\& + \frac{{41\gamma _{ij}^3}}{{240{h_i}{h_j}}} - \frac{{41{\gamma _{ij}}{\partial _i}{\gamma _{ij}}}}{{1440h_i^2}} + \frac{{{\partial _i}{\gamma _{ij}}{\partial _j}{h_j}}}{{1440{h_i}h_j^2}} + \frac{{79\gamma _{ij}^2{\partial _j}{h_j}}}{{1440h_j^3}} - \frac{{{\gamma _{ij}}{\partial _i}{h_i}{\partial _j}{h_j}}}{{720h_i^2h_j^2}} - \frac{{{\gamma _{ij}}{h_i}{{({\partial _j}{h_j})}^2}}}{{288h_j^5}}
\\& + \sum\limits_k {} (\frac{{{\gamma _{ij}}{\gamma _{ik}}{H_j}}}{{60{h_j}{h_k}}} - \frac{{{\gamma _{ik}}{\gamma _{jk}}{h_i}{h_j}{H_k}}}{{30h_k^4}} - \frac{{{\gamma _{ij}}{\gamma _{jk}}{h_i}{H_j}}}{{60h_j^2{h_k}}} + \frac{{{\gamma _{ik}}{\gamma _{jk}}{h_i}{H_j}}}{{60{h_j}h_k^2}} - \frac{{7{\gamma _{ij}}{\gamma _{jk}}{h_i}{H_k}}}{{60h_j^2{h_k}}}
\\& - \frac{{{\gamma _{ij}}{\gamma _{ik}}{\partial _j}{h_j}}}{{720h_j^2{h_k}}} + \frac{{{\gamma _{ij}}{\gamma _{jk}}{h_i}{\partial _j}{h_j}}}{{240h_j^3{h_k}}} - \frac{{{\gamma _{ik}}{\gamma _{jk}}{h_i}{\partial _j}{h_j}}}{{1440h_j^2h_k^2}} + \frac{{{\gamma _{ij}}{\gamma _{jk}}{h_i}{\partial _k}{h_k}}}{{720h_k^4}} + \frac{{{\gamma _{ik}}{\gamma _{jk}}{h_i}{h_j}{\partial _k}{h_k}}}{{288h_k^5}}
\\& + \frac{{{\gamma _{jk}}{\partial _i}{\gamma _{ij}}}}{{1440{h_i}{h_k}}} - \frac{{{h_j}{h_k}{\gamma _{ij}}{\partial _i}{\gamma _{ik}}}}{{360h_i^4}} - \frac{{{h_j}(3{\gamma _{ik}}{\partial _i}{\gamma _{ij}} + 2{\gamma _{ij}}{\partial _i}{\gamma _{ik}})}}{{1440h_i^2{h_k}}} - \frac{{7{h_j}{\gamma _{ij}}{\partial _k}(h_k^{ - 1}{\gamma _{ik}})}}{{1440h_i^2}}
\\& - \frac{{{h_i}{h_j}{\gamma _{ik}}{\partial _k}{\gamma _{jk}}}}{{480h_k^4}} + \frac{{\gamma _{ij}^2{\gamma _{jk}}}}{{120{h_j}{h_k}}} + \frac{{7{h_i}{\gamma _{ij}}\gamma _{jk}^2}}{{160h_j^3}} + \frac{{11{\gamma _{ij}}{\gamma _{ik}}{\gamma _{jk}}}}{{2880h_k^2}} + \frac{{{h_j}\gamma _{ik}^2{\gamma _{jk}}}}{{96h_k^3}})
\\& + \sum\limits_{k,l} {} (\frac{{{h_i}{h_j}{\gamma _{il}}{\gamma _{jl}}}}{{720{h_k}h_l^2}}(\frac{{{\gamma _{kl}}}}{{{h_l}}} - \frac{{{\gamma _{jk}}}}{{2{h_j}}}) - \frac{{{h_i}{\gamma _{ij}}{\gamma _{jl}}{\gamma _{kl}}}}{{720{h_k}h_l^2}}),
  \end{align*}
and
 \begin{align*}
Q_i^{(2)} =& \frac{{4H_i^3}}{{5h_i^2}} - \frac{{7{\partial _i}{h_i}H_i^2}}{{10h_i^3}} + \frac{{7{{({\partial _i}{h_i})}^2}{H_i}}}{{48h_i^4}} - \frac{{{{({\partial _i}{h_i})}^3}}}{{120h_i^5}} + \sum\limits_k {} (\frac{{7{\gamma _{ik}}{H_i}{H_k}}}{{10{h_i}{h_k}}} - \frac{{{\gamma _{ik}}{\partial _i}{h_i}{H_i}}}{{120h_i^2{h_k}}}
\\& + \frac{{7{\partial _k}(h_k^{ - 1}{\gamma _{ik}}){H_i}}}{{240{h_i}}} - \frac{{7{\gamma _{ik}}{\partial _i}{h_i}{H_k}}}{{80h_i^2{h_k}}} + \frac{{{\gamma _{ik}}{H_k}}}{{576{u_{ik}}{h_i}{h_k}}} + \frac{{(2{H_i} + 7{H_k}){\partial _i}{\gamma _{ik}}}}{{240{h_i}{h_k}}}
\\& + \frac{{{\gamma _{ik}}{h_k}{H_i}}}{{576{u_{ik}}h_i^3}} - \frac{{31\gamma _{ik}^2{H_i}}}{{144h_i^2}} + \frac{{{\gamma _{ik}}{{({\partial _i}{h_i})}^2}}}{{720h_i^3{h_k}}} + \frac{{253\gamma _{ik}^2{\partial _i}{h_i}}}{{5760h_i^3}} - \frac{{{\partial _i}{\gamma _{ik}}{\partial _i}{h_i}}}{{960h_i^2{h_k}}} - \frac{{\gamma _{ik}^2{\partial _k}{h_k}}}{{2880h_k^3}}
\\& - \frac{{7{\partial _k}(h_k^{ - 1}{\gamma _{ik}}){\partial _i}{h_i}}}{{1920h_i^2}} - \frac{{7{\partial _i}{\gamma _{ik}}{\partial _k}{h_k}}}{{5760{h_i}h_k^2}} - \frac{{41{\partial _i}{\gamma _{ik}}{\partial _i}{h_i}{h_k}}}{{5760h_i^4}} + \frac{{{\partial _i}({h_i}{\gamma _{ik}}){\partial _k}{h_k}}}{{2880h_k^4}}
\\& - \frac{{113{\gamma _{ik}}{\partial _i}{\gamma _{ik}}}}{{5760h_i^2}} + \frac{{(3{\partial _i}{\gamma _{ik}} + {\partial _k}{\gamma _{ik}}){\gamma _{ik}}}}{{1440h_k^2}} - \frac{{{\partial _i}{\gamma _{ik}}{h_k}}}{{576{u_{ik}}h_i^3}} - \frac{{{\partial _k}{\gamma _{ik}}}}{{576{u_{ik}}{h_i}{h_k}}} - \frac{{\gamma _{ik}^3}}{{240{h_i}{h_k}}})
\\& + \sum\limits_{k,l} {} ( - \frac{{{\gamma _{kl}}{\partial _i}({h_i}{\gamma _{il}})}}{{2880{h_k}h_l^2}} + \frac{{\gamma _{il}^2{\gamma _{kl}}}}{{2880{h_k}{h_l}}} - \frac{{{\gamma _{ik}}\gamma _{il}^2}}{{240{h_i}{h_k}}} - \frac{{{\gamma _{kl}}{\partial _i}{\gamma _{ik}}}}{{2880{h_i}{h_l}}} + \frac{{{u_{lk}}{\gamma _{ik}}{\partial _l}{\gamma _{kl}}}}{{1152{u_{il}}{h_i}{h_l}}}
\\& + \frac{{{u_{kl}}{\gamma _{ik}}{\gamma _{kl}}{\partial _i}{\gamma _{il}}}}{{144h_i^2}} + \frac{{{h_l}{\gamma _{ik}}{\partial _i}{\gamma _{il}}}}{{144h_i^2{h_k}}} + \frac{{{h_k}{u_{kl}}{\gamma _{kl}}{\partial _i}{\gamma _{il}}}}{{1152{u_{ik}}h_i^3}} + \frac{{{h_l}{u_{ik}}\gamma _{ik}^2{\partial _i}{\gamma _{il}}}}{{40h_i^3}}).
                          \end{align*}

In these expressions, all summations are taken over the ranges of indices where the denominators do not vanish.


\vspace{30pt} \noindent
Xiaobo Liu \\
Beijing International Center for Mathematical Research,\\
Beijing University, Beijing, China.\\
E-mail address: {\it xbliu@math.pku.edu.cn} \\
\& \\
\noindent
Department of Mathematics, \\
University of Notre Dame, \\
Notre Dame,  IN  46556, USA \\
E-mail address: {\it xliu3@nd.edu}\\
\\
\\
Xin Wang \\
School of Mathematical Sciences,  \\
Beijing University, Beijing, China \\
E-mail address:{\it xinwang-1989@163.com}


\begin{thebibliography}{399}
\bibitem[AGV]{AGV}
D. Abramovich, T. Graber, A. Vistoli, {\it Gromov-Witten theory of
Deligne-Mumford stacks}, Amer. J. Math., Volume 130, Number 5,
(2008), pp. 1337-1398.

\bibitem[CR]{CR} W.Chen and Y.Ruan,
{\it Orbifold Gromov-Witten theory}, Orbifolds in mathematics and physics. Contemp. Math., 310,
Amer. Math. Soc., Providence, RI(2002):25-85.

\bibitem[D]{D}   B. Dubrovin,
        {\it Geometry of 2D topological field theories}, Integrable Systems and Quantum Group,
        Spinger Lecture Notes in Math. 1620 (1996), 120-348.

\bibitem[DLZ]{DLZ}   B. Dubrovin, S. Liu, Y. Zhang,
        {\it On the genus Two Free Energies for Semisimple Frobenius Manifold},
        Russian Journal of Mathematical Physics 19 (2012), 273-298.
\bibitem[FLZZ]{FLZZ}   Y. Fu, S. Liu, Y. Zhang, C. Zhou,
        {\it Proof of a Conjecture On the genus Two Free Energies Associated to
        the $A_{n}$ Singularity},
        arXiv:1305.1008.


\bibitem[KM]{KM} M. Konsevich, Yu. I. Manin,
        {\it Gromov-Witten classes, quantum cohomology, and enumerative geometry},
       Commun. Math. Phys. 164 (1994), 525-562.

\bibitem[KS]{KS}  M. Krawitz, Y. Shen, {\it Landau-Ginzburg/Calabi-Yau Correspondence of all Genera for Elliptic Orbifold $\mathbb{P}^{1}$}, arXiv:1106.6270

\bibitem[L02]{L02} X. Liu,
    {\it Quantum product on the big phase space and Virasoro conjecture},
    Advances in Mathematics 169 (2002), 313-375.

\bibitem[L06]{L06} X. Liu,
    {\it Idempotents on the big phase space},
    Contemporary Mathematics,vol. 403 (2006), 43-66.

\bibitem[L07]{L07} X. Liu,
    {\it Genus-2 Gromov-Witten invariants for manifolds with semisimple quantum cohomology},
    American Journal of Mathematics, 129 (2007), no.2, 463-498.
\bibitem[LW]{LW}   X.Liu, X.Wang,
        {\it Conditions for the vanishing of the genus-2 G-function},
        arXiv:math/1310.2101.

\end{thebibliography}
\end{document}